\documentclass[a4paper,11pt]{article}
\usepackage{amsmath,amssymb,amsthm}
\usepackage[no-math]{fontspec}
\usepackage{xeCJK}
\setCJKsansfont{ipaexg.ttf}[
  Scale=0.962216,
  BoldFont=ipaexg.ttf,
  ItalicFont=ipaexg.ttf,
  BoldItalicFont=ipaexg.ttf
]
\setCJKmonofont{ipaexg.ttf}[
  Scale=0.962216,
  BoldFont=ipaexg.ttf
]
\usepackage{xcolor}
\definecolor{DOIblue}{RGB}{0,80,160}
\usepackage[unicode=true,colorlinks=true,linkcolor=black,citecolor=black,urlcolor=DOIblue]{hyperref}
\hypersetup{
  pdftitle={Uniform counting of quadratic fields with local conditions},
  pdfauthor={Takashi Taniguchi},
  pdfsubject={English translation followed by the Japanese original},
  bookmarksnumbered=true
}
\AtBeginDocument{\lefthyphenmin=2\relax\righthyphenmin=3\relax}

\normalfont

\makeatletter
\newcommand{\bilingualmaketitle}{%
  \par\begingroup
    \renewcommand\thefootnote{\@fnsymbol\c@footnote}%
    \def\@makefnmark{\rlap{\@textsuperscript{\normalfont\@thefnmark}}}%
    \long\def\@makefntext##1{\parindent 1em\noindent
      \hb@xt@1.8em{\hss\@textsuperscript{\normalfont\@thefnmark}}##1}%
    \newpage\global\@topnum\z@
    \@maketitle
    \thispagestyle{plain}\@thanks
  \endgroup
  \setcounter{footnote}{0}%
  \global\let\@thanks\@empty
}
\makeatother

\newcommand{\bilingualtheoremname}{Theorem}
\newcommand{\bilinguallemmaname}{Lemma}
\newcommand{\bilingualpropositionname}{Proposition}
\newcommand{\bilingualcorollaryname}{Corollary}
\newcommand{\bilingualremarkname}{Remark}
\newtheorem{theorem}{\bilingualtheoremname}
\newtheorem{lemma}[theorem]{\bilinguallemmaname}
\newtheorem{proposition}[theorem]{\bilingualpropositionname}
\newtheorem{corollary}[theorem]{\bilingualcorollaryname}
\newtheorem{remark}[theorem]{\bilingualremarkname}

\newcommand{\bilingualversion}{en}

\newcommand{\Q}{\mathbb{Q}}
\newcommand{\Z}{\mathbb{Z}}
\newcommand{\leg}[2]{\left(\frac{#1}{#2}\right)}
\newcommand{\Disc}{\operatorname{Disc}}

\newcommand{\address}[1]{%
  \par\medskip\noindent{\small #1}\par}
\newcommand{\email}[1]{%
  \noindent{\small\href{mailto:#1}{\texttt{#1}}}\par}

\title{Uniform counting of quadratic fields with local conditions}
\author{Takashi Taniguchi}
\date{}

\begin{document}
\pdfbookmark[0]{English translation}{en:version}
\bilingualmaketitle


\begin{abstract}
\setlength{\emergencystretch}{1em}
We sharpen the uniform counting estimates of Bhargava, Thorne,
and the author for quadratic fields with prescribed local conditions.
Combining the Burgess and P\'olya--Vinogradov bounds, we improve the
dependence of the error term on the primes supporting the local conditions
and make explicit the factors previously absorbed into an arbitrary
$\epsilon$-power of their product.
The resulting estimate has an absolute implied constant.
\end{abstract}

\section{Introduction}\label{en:sec:statement}

For each prime $p$, let $A_p$ denote the set of isomorphism classes of
\'etale quadratic $\Q_p$-algebras.
Let $n$ be a positive squarefree integer.
A family of nonempty proper subsets
$\sigma_p$ of $A_p$, one for each prime $p\mid n$,
is called a \emph{local specification} for quadratic fields
and is denoted by $\sigma=(\sigma_p)_{p\mid n}$.
We call $n$ the \emph{support} of $\sigma$.
A quadratic field $F$ is said to satisfy $\sigma$ if
$F\otimes \Q_p\in\sigma_p$ for every $p\mid n$.
Write $\mathcal F(\sigma)$ for the set of isomorphism classes of
quadratic fields satisfying $\sigma$, and, for each sign $\pm$ and $X>0$, put
\[
 N_2^\pm(X,\sigma)
 :=\#\{F\in \mathcal F(\sigma)\mid 0<\pm\Disc(F)<X\}.
\]
Our main theorem gives an estimate for this counting function that is
uniform in $X$ and $\sigma$.
We call $\sigma$ \emph{pure} if, for every odd prime $p\mid n$,
$\sigma_p$ does not contain both an unramified algebra
and a ramified algebra.
In this section, we state the main theorem for pure local specifications,
which are particularly important in applications.
For general local specifications, see Corollary~\ref{en:cor:general}
in Section~\ref{en:sec:local}.

To state the theorem, we define several quantities associated with a
pure local specification $\sigma$.
For $i=1,2$, let $u_i$ be the product of the odd primes $p\mid n$
for which $\sigma_p$ consists of $i$ unramified algebras, and let
$r_i$ be the product of the odd primes $p\mid n$ for which $\sigma_p$
consists of $i$ ramified algebras.
Thus $n=u_1u_2r_1r_2$ if $n$ is odd, and
$n=2u_1u_2r_1r_2$ if $n$ is even.
For each odd prime $p\mid n$, define the local density $c_p(\sigma_p)$ by
$c_p(\sigma_p):=ip/(2(p+1))$ if $p\mid u_i$, and by
$c_p(\sigma_p):=i/(2(p+1))$ if $p\mid r_i$ ($i=1,2$).
If $2\mid n$, assign the local density $1/(3\cdot 2^e)$ to each
\'etale quadratic $\Q_2$-algebra whose discriminant has $2$-adic
valuation $e$, and define $c_2(\sigma_2)$ as the sum of these
densities over the elements of $\sigma_2$.

For a positive squarefree integer $q$, let $\omega(q)$ denote the number
of distinct prime divisors of $q$, and put
$\lambda_2(q):=\prod_{p\mid q}(1+p^{-1/2})$.
Our main result is as follows.
\begin{theorem}\label{en:thm:main}
For every pure local specification $\sigma$ and $X>0$, we have
\begin{equation}\label{en:eq:main}
 N_2^\pm(X,\sigma)
 =\frac{3}{\pi^2}X\prod_{p\mid n}c_p(\sigma_p)
 +O\!\left(
 \frac{X^{1/2}\,u_1^{3/16}}{r_1^{5/16}r_2^{1/2}}\,
 \left(\frac32\right)^{\omega(u_1r_1)}
 (\log(2u_1r_1))^{5/4}\lambda_2(u_2r_2)
 \right).
\end{equation}
The implied constant is absolute.
\end{theorem}

When $n=1$, so that no local conditions are imposed, the count reduces
to Gegenbauer's count of squarefree integers~\cite{en:Gegenbauer}.
For a fixed local specification $\sigma$, one can count squarefree
integers in residue classes as in~\cite{en:Gegenbauer}, obtaining
\eqref{en:eq:main} with an error term $O_\sigma(X^{1/2})$.
Our aim is to obtain an error bound with more precise dependence
on $\sigma$.

The dependence of counting estimates on local conditions has also
been studied extensively.
For counts with coprimality conditions, corresponding to the case
$u_1=r_1=r_2=1$ here, Cohen~\cite[Lemma 5.2]{en:Cohen} gave a uniform
estimate, and Suryanarayana~\cite[Corollary 1]{en:Suryanarayana} and
Robert and Tenenbaum~\cite[(10.1)]{en:RT} refined the dependence of the
error term on the modulus.
Theorem~\ref{en:thm:main} further improves that dependence in this case.
For counting squarefree integers in residue classes,
Cohen and Robinson~\cite[Corollary 1]{en:CR} and
Schwarz~\cite[Lemma 8]{en:Schwarz} obtained estimates uniform in the
modulus and the residue class.
Bhargava, Thorne, and the author~\cite[Proposition 4.2]{en:BTT}
treated arbitrary $\sigma$, obtaining, for pure $\sigma$, the bound
\begin{equation}\label{en:eq:BTT}
 O_{\epsilon}\!\left(
 X^{1/2}\,u_1^{1/4}u_2^{1/4}r_1^{-1/4}r_2^{-1/2}
 (u_1u_2r_1r_2)^{\epsilon}
 \right)
\end{equation}
for the error term in \eqref{en:eq:main}, where $\epsilon>0$ is arbitrary
and the implied constant depends only on $\epsilon$.
This estimate permits the specification of individual ramified local
extensions while retaining the saving $r_1^{-1/4}r_2^{-1/2}$
from the ramified primes.
Our estimate \eqref{en:eq:main} improves the power factors
$u_1^{1/4}u_2^{1/4}r_1^{-1/4}$ to $u_1^{3/16}r_1^{-5/16}$
and makes explicit the factors absorbed into the $\epsilon$-power
of the support in \eqref{en:eq:BTT}.

The key to the proof of Theorem~\ref{en:thm:main} is an estimate for the
character sums over squarefree integers in \eqref{en:eq:T_h}.
Expressing the indicator of the local condition at each prime
$p\mid u_1r_1$ in terms of quadratic characters and expanding the
product leads to such sums.
For sums of this type, Granville and Soundararajan~\cite[Lemma 4.1]{en:GS}
obtained a bound using the P\'olya--Vinogradov inequality, while
Munsch~\cite[Lemma 2.3]{en:Munsch} obtained one using
Burgess's inequality~\cite{en:Burgess}.
After M\"obius inversion, we estimate the inner character sums in
\eqref{en:eq:T_h-Mobius} by the P\'olya--Vinogradov inequality or the
Burgess bound of Jain-Sharma, Khale, and Liu~\cite[Theorem 1.1]{en:JKL},
according to the length of the summation interval.
We track the dependence on the local conditions through explicit
prime-by-prime factors.

Uniform estimates for counting cubic fields with local conditions,
due to Thorne and the author~\cite{en:TT} and
Bhargava, Thorne, and the author~\cite{en:BTT}, have found applications to
various problems in arithmetic statistics.
Such arguments sometimes require counts of lower-degree fields with
local conditions to isolate the contributions of reducible \'etale
algebras.
For example, in~\cite[\S4]{en:BTT}, counts of quadratic fields
with local conditions are used to estimate the contribution
of reducible cubic \'etale algebras $\Q\times F$,
where $F/\Q$ is quadratic.
Shankar and Tsimerman~\cite[\S3.1]{en:ST} likewise use such counts to
estimate the contribution of reducible quartic \'etale algebras
in counting quartic fields.
Our result sharpens the uniform estimates for the quadratic-field
counts used in these arguments.
Other applications of counts of quadratic fields with local conditions
include Pollack's~\cite[\S5.2]{en:Pollack} study of the average least
inert prime in quadratic fields and bounds on $\ell$-torsion in
class groups obtained by Ellenberg, Pierce, and Wood~\cite[\S6]{en:EPW}.
In the latter work, the quadratic field counts serve as input to a sieve.
We expect Theorem~\ref{en:thm:main} and its generalization,
Corollary~\ref{en:cor:general}, to find further applications in
arithmetic statistics.

The paper is organized as follows.
In Section~2, we combine the Burgess bound and the
P\'olya--Vinogradov inequality to obtain a uniform bound for a sum
of absolute values of character sums.
In Section~3, we use this bound and M\"obius inversion to derive
a uniform counting formula for squarefree integers with prescribed
congruence conditions and Legendre-symbol values.
In Section~4, we reduce the counting of quadratic-field discriminants
to the problem studied in Section~3 and prove Theorem~\ref{en:thm:main}.
We also deduce Corollary~\ref{en:cor:general} for general local specifications.
Unless otherwise stated, the implied constants in Landau's $O$-notation
and Vinogradov's $\ll$-notation are absolute.

\subsection*{Acknowledgments}
This work was supported by JSPS KAKENHI Grant Numbers
JP23K22386 and JP24K21512, and by the Research Institute
for Mathematical Sciences (RIMS), an International Joint
Usage/Research Center at Kyoto University.
Part of this work was carried out during the author's
participation in the Thematic Programme on Arithmetic
Statistics at the Lodha Mathematical Sciences Institute.
The author is deeply grateful to the organizers and staff
of the programme for providing a truly exceptional
research environment.

\subsection*{Note on translation and AI assistance}
This paper was originally written in Japanese, the author's native language. The author subsequently translated it into English with the assistance of ChatGPT (GPT-6) in October and November 2026. During this process, the author also consulted ChatGPT on mathematical notation and terminology.
The author has carefully reviewed and revised the entire translation and takes full responsibility for the content of the paper.
Both the Japanese and English versions are available on the author's website.

\section{Character sum estimates}\label{en:sec:characters}

Throughout this section, $\chi$ denotes
a nontrivial primitive Dirichlet character of conductor $q$.
For a nonnegative integer $N$ and a real number $b$, define the
$N$-term character sum
\[
 S(\chi,N;b):=\sum_{b<j\le b+N}\chi(j).
\]
For a positive integer $q$, write $\omega(q)$ for the number of
distinct prime factors of $q$, $\tau(q)$ for the number of positive
divisors of $q$, and $\varphi(q)$ for the order of
$(\Z/q\Z)^\times$. Define
\begin{equation}\label{en:eq:kappa}
 \kappa(q):=q^{3/16}(2^{\omega(q)}\tau(q))^{3/4}
 \left(\frac q{\varphi(q)}\right)^{1/2}.
\end{equation}
Since $2^{\omega},\tau,\varphi$ are multiplicative, so is $\kappa$.
We use the following explicit form of Burgess's bound~\cite{en:Burgess}
for character sums, due to
Jain-Sharma, Khale, and Liu~\cite[Theorem 1.1]{en:JKL}.
\begin{lemma}[Burgess bound \cite{en:JKL}]\label{en:lem:all-q}
For all $N\le q^{5/8}$, we have
\begin{equation}\label{en:eq:JKL-usable}
 |S(\chi,N;b)|\ll\sqrt N\,\kappa(q)(\log q)^{1/4}.
\end{equation}
\end{lemma}
The assumption in~\cite[Theorem 1.1]{en:JKL} that the conductor is
sufficiently large can be removed by enlarging the implied constant.
The next lemma is the familiar P\'olya--Vinogradov inequality;
see, for example, Davenport~\cite[\S23]{en:Davenport} for a proof.
\begin{lemma}[P\'olya--Vinogradov]\label{en:lem:PV}
For all $\chi,N,b$, we have
\[
 |S(\chi,N;b)|\ll\sqrt q\,\log q.
\]
\end{lemma}

Combining the Burgess and P\'olya--Vinogradov bounds,
we obtain the following proposition.
\begin{proposition}\label{en:prop:aggregate}
Let $Y>0$. Suppose that, for each integer $a\ge1$, $N_a$ is a
nonnegative integer satisfying $N_a\le Y/a^2$. Then
\begin{equation}\label{en:eq:aggregate}
 \sum_{a\ge1}|S(\chi,N_a;b_a)|
 \ll\sqrt Y\,\kappa(q)(\log q)^{5/4},
\end{equation}
where the $b_a$ are arbitrary real numbers.
\end{proposition}
\begin{proof}
Put
$
 A=\sqrt Y\,q^{-5/16}.
$
For $a\leq A$, Lemma~\ref{en:lem:PV} gives
\[
 \sum_{a\leq A}|S(\chi,N_a;b_a)|
 \ll A\sqrt q\,\log q
 =\sqrt Y\,q^{3/16}\log q
 \ll\sqrt Y\,\kappa(q)(\log q)^{5/4}.
\]
If $a>A$, then $N_a\le Y/A^2=q^{5/8}$.
If $N_a>0$, then $N_a\ge1$ by integrality, and hence
$a\le\sqrt Y$.
Thus Lemma~\ref{en:lem:all-q} gives
\[
 \sum_{a> A}|S(\chi,N_a;b_a)|
 \ll     \sum_{A< a\le\sqrt Y}
 \sqrt{N_a}\,\kappa(q)(\log q)^{1/4}
 \leq\sqrt Y\,\kappa(q)(\log q)^{1/4}
      \sum_{A< a\le\sqrt Y}\frac1a.
\]
The last sum is bounded by
\[
 \sum_{A< a\le\sqrt Y}\frac1a
 \le1+\log\frac{\sqrt Y}{A}
 =1+\frac5{16}\log q.
\]
The bounds for $a\le A$ and $a>A$ together prove the proposition.
\end{proof}

\section{Counting squarefree integers}\label{en:sec:sieve}

Let $n_1,n_2$ be coprime positive odd squarefree integers, and let
$c\in(\Z/8\Z)^\times$.
For each $p\mid n_1$, choose a sign $\varepsilon(p)\in\{+1,-1\}$,
and write $\varepsilon=(\varepsilon(p))_{p\mid n_1}$.
For $Y>0$, let $T(Y;c,n_1,n_2,\varepsilon)$ be the number of integers
$0<d<Y$ satisfying the following conditions:
\begin{center}
$d$ is squarefree;
\quad
$d\equiv c\pmod8$;
\quad
$(d,n_1n_2)=1$;
\quad
$\leg dp=\varepsilon(p)$ for every $p\mid n_1$.
\end{center}
By definition, this counting function is
\begin{equation}\label{en:eq:T-definition}
 T(Y;c,n_1,n_2,\varepsilon)=
 \sum_{\substack{0<d<Y\\ d\equiv c\ (8)\\ (d,n_1n_2)=1}}
 \mu(d)^2\prod_{p\mid n_1}\frac{1+\varepsilon(p)\leg dp}{2}.
\end{equation}
Here $\mu$ is the M\"obius function.
Define a multiplicative function $\lambda_1$ on the positive squarefree
integers by prescribing, for each prime $p$,
\begin{equation}\label{en:eq:lambda-def}
 \lambda_1(p):=\frac{\kappa(p)+\lambda_2(p)}2
 =\frac{\sqrt2}{(1-p^{-1})^{1/2}}\,p^{3/16}+\frac{1+p^{-1/2}}2.
\end{equation}
We prove the following uniform estimate.
\begin{proposition}\label{en:prop:sieve}
\begin{equation}\label{en:eq:sieve}
 T(Y;c,n_1,n_2,\varepsilon)=\frac{Y}{2^{\omega(n_1)}\pi^2}\prod_{p\mid n_1n_2}\frac p{p+1}
 +O\!\left(\sqrt Y\,(\log(2n_1))^{5/4}\lambda_1(n_1)\lambda_2(n_2)\right).
\end{equation}
\end{proposition}

\begin{proof}
Extend $\varepsilon$ to a multiplicative function on the divisors of
$n_1$ by putting
$\varepsilon(h):=\prod_{p\mid h}\varepsilon(p)$ for each $h\mid n_1$.
The product in \eqref{en:eq:T-definition} expands as
\[
 \prod_{p\mid n_1}\frac{1+\varepsilon(p)\leg dp}{2}
 =\frac{1}{2^{\omega(n_1)}}\sum_{h\mid n_1}\varepsilon(h)\leg dh.
\]
Here $\leg dh$ is the Jacobi symbol.
For each $h\mid n_1$, put
\begin{equation}\label{en:eq:T_h}
 T_h(Y):=
 \sum_{\substack{0<d<Y\\d\equiv c\ (8)\\ (d,n_1n_2)=1}}
 \mu(d)^2\leg dh.
\end{equation}
Then
\begin{equation}\label{en:eq:sieve-expand}
 T(Y;c,n_1,n_2,\varepsilon)
 =\frac{1}{2^{\omega(n_1)}}\sum_{h\mid n_1}\varepsilon(h)T_h(Y).
\end{equation}
We suppress the dependence of $T_h(Y)$ on $c,n_1,n_2$ from the notation.
We now substitute $\mu(d)^2=\sum_{a^2\mid d}\mu(a)$ into
$T_h(Y)$ and put $d=a^2k$.
The original conditions give
$(a,2n_1n_2)=(k,n_1n_2)=1$; since
$a^2\equiv1\pmod8$ and $\leg{a^2}h=1$, we obtain
\[
 T_h(Y)=
 \sum_{\substack{a\le\sqrt Y\\ (a,2n_1n_2)=1}}\mu(a)
   \sum_{\substack{0<k<Y/a^2\\ k\equiv c\ (8)\\ (k,n_1n_2)=1}}\leg kh.
\]
As $\leg kh=0$ if $(k,h)>1$, the condition $(k,n_1n_2)=1$ in
the inner sum may be replaced by $(k,n_1n_2/h)=1$.
Inclusion--exclusion gives
\[
 \mathbf 1_{(k,n_1n_2/h)=1}
 =\sum_{e\mid (n_1n_2/h,k)}\mu(e)
 =\sum_{\substack{e\mid n_1n_2/h\\e\mid k}}\mu(e).
\]
Put $k=e\ell$. Since $e$ is odd, it has an inverse $e^{-1}$ in
$(\Z/8\Z)^\times$, and $k\equiv c\pmod8$ becomes
$\ell\equiv ce^{-1}\pmod8$.
Interchanging the sums over $e$ and $a$ gives
\begin{equation}\label{en:eq:T_h-Mobius}
 T_h(Y)=\sum_{e\mid n_1n_2/h}\mu(e)\leg eh
   \sum_{\substack{a\le\sqrt{Y/e}\\ (a,2n_1n_2)=1}}\mu(a)
   \sum_{\substack{0<\ell<Y/ea^2\\\ell\equiv ce^{-1}\ (8)}}\leg\ell h.
\end{equation}

We first estimate $T_1(Y)$.
In this case, the sum over $\ell$ in \eqref{en:eq:T_h-Mobius} is
$Y/(8ea^2)+O(1)$.
On summing over $a$, both the total contribution of the $O(1)$
terms and the error from extending the sum over $a$ in the main term
to infinity are $O(\sqrt{Y/e})$.
Thus, since $\sum_{e\mid n_1n_2}e^{-1/2}=\lambda_2(n_1n_2)$, we have
\begin{align*}
T_1(Y)&=\frac{Y}{8}
   \sum_{e\mid n_1n_2}\frac{\mu(e)}e
   \sum_{\substack{a\ge1\\ (a,2n_1n_2)=1}}\frac{\mu(a)}{a^2}
   +O\!\left(\sqrt Y\,\lambda_2(n_1n_2)\right)
 \\
 &=
 \frac{Y}{8}
 \prod_{p\mid n_1n_2}(1-p^{-1})\prod_{p\nmid2n_1n_2}(1-p^{-2})
 +O\!\left(\sqrt Y\,\lambda_2(n_1n_2)\right).
\end{align*}
Since $\prod_{p\ne2}(1-p^{-2})=8/\pi^2$, the main term of
$2^{-\omega(n_1)}T_1(Y)$ agrees with the main term of \eqref{en:eq:sieve}.

We next estimate $T_h(Y)$ for $h>1$.
For each $e\mid n_1n_2/h$, choose an integer $r$ satisfying
$r\equiv ce^{-1}\pmod8$, and an integer $b$ satisfying
$8b\equiv r\pmod h$.
Let $N_a$ denote the number of integers $\ell$ such that
$0<\ell<Y/ea^2$ and $\ell\equiv r\pmod8$.
Writing $\ell=r+8j$, we have $\ell>0$ if and only if $j>-r/8$,
and $\leg\ell h=\leg8h\leg{j+b}h$.
Consequently,
\[
 \sum_{\substack{0<\ell<Y/ea^2\\\ell\equiv ce^{-1}\ (8)}}\leg\ell h
 =\leg8h\,S\!\left(\leg{\cdot}h,N_a;b-\frac r8\right),
 \qquad N_a\le\frac{Y}{ea^2}.
\]
Since the Jacobi symbol $\leg \cdot h$ is a primitive character
of conductor $h$, Proposition~\ref{en:prop:aggregate} gives
\[
 \sum_{\substack{a\le\sqrt{Y/e}\\ (a,2n_1n_2)=1}}
 \left|S\!\left(\leg{\cdot}h,N_a;b-\frac r8\right)\right|
 \ll\sqrt{\frac Ye}\,\kappa(h)(\log h)^{5/4}.
\]
Using $h<2n_1$ and
$\sum_{e\mid n_1n_2/h}e^{-1/2}=\lambda_2(n_1n_2/h)$, we obtain
\begin{equation}\label{en:eq:nonprincipal}
\sum_{\substack{1\neq h\mid n_1}}|T_h(Y)|
 \ll \sqrt Y\,(\log(2n_1))^{5/4}
     \sum_{\substack{1\neq h\mid n_1}}\kappa(h)\lambda_2(n_1n_2/h).
\end{equation}
Since $\kappa(1)=1$, the error bound for $T_1(Y)$ and
\eqref{en:eq:nonprincipal} can be combined into a sum over all $h\mid n_1$.
By the multiplicativity of $\lambda_2$ and $\kappa$
and the definition \eqref{en:eq:lambda-def}, we have
\begin{equation}\label{en:eq:sieve-product}
 \sum_{h\mid n_1}\kappa(h)\lambda_2(n_1n_2/h)
 =\lambda_2(n_2)
 \prod_{p\mid n_1}\bigl(\lambda_2(p)+\kappa(p)\bigr)
 =2^{\omega(n_1)}\lambda_1(n_1)\lambda_2(n_2).
\end{equation}
Together with \eqref{en:eq:sieve-expand}, this proves
\eqref{en:eq:sieve}.
\end{proof}

\section{Proof of the main theorem and general local specifications}\label{en:sec:local}
We first prove Theorem~\ref{en:thm:main} for a pure local specification
$\sigma$, and then deduce Corollary~\ref{en:cor:general} for general $\sigma$.
\begin{proof}[Proof of Theorem~\ref{en:thm:main}]
We first describe the local condition at $2$ in terms of the
discriminant $D=\Disc(F)$ of a quadratic field $F$.
Put
\[
 e=\operatorname{ord}_2(D)\in\{0,2,3\},\qquad
 c\equiv D/2^e\pmod8,\quad c\in(\Z/8\Z)^\times.
\]
Then the isomorphism class of $F\otimes\Q_2$ is determined by $(e,c)$.
There are exactly eight possible pairs $(e,c)$:
\[
(0,1),(0,5),(2,3),(2,7),(3,1),(3,3),(3,5),(3,7).
\]
The density of the corresponding quadratic $\Q_2$-algebra is
$1/(3\cdot2^e)$ in each case.

It suffices to prove the theorem when $2\mid n$ and
$\sigma_2$ is a singleton.
Indeed, partitioning according to the isomorphism class of
$F\otimes\Q_2$ (and replacing $n$ by $2n$ if necessary)
expresses the general count as a sum of at most eight
counts of this form, without changing $u_1,u_2,r_1,r_2$.
The main terms add by additivity of local densities
(which sum to $1$ over $A_2$), and the error bound
is preserved up to an absolute factor.
Write $(e,c)$ for the pair corresponding to the unique
element of $\sigma_2$.
Then $c_2(\sigma_2)=1/(3\cdot 2^e)$.

The product of the odd primes at which ramification is prescribed is
$r_1r_2$.
Thus the discriminant $D$ of any quadratic field $F$ counted by
$N_2^\pm(X,\sigma)$ can be written as
$D=\pm 2^e r_1r_2d$ for some integer $d$.
We express the conditions on $D$ in terms of $d$.
Putting $Y=X/(2^e r_1r_2)$, we have
\[
\text{$d$ is squarefree};
\qquad
0<d<Y;
\qquad
(d,n)=1;
\qquad
d\equiv\pm c(r_1r_2)^{-1}\pmod8.
\]
Among the odd primes $p\mid n$, the product of those with
$\#\sigma_p=1$ is $u_1r_1$, and the product of those with
$\#\sigma_p=2$ is $u_2r_2$.
For $p\mid u_1r_1$, we also have a condition
$\leg dp=\varepsilon(p)$.
Indeed, $F\otimes\Q_p$ is determined by $\leg Dp$ if
$p\mid u_1$, and by $\leg{D/p}p$ if $p\mid r_1$.
Substituting $D=\pm 2^e r_1r_2d$ gives a condition of the form
$\leg dp=\pm1$ in either case.
For $p\mid u_2r_2$, no condition remains beyond the already imposed
$(d,p)=1$.
Conversely, the integers $d$ satisfying these conditions correspond
bijectively to the quadratic fields $F$ being counted, apart from the
case $D=1$.
Hence
\[
N_2^\pm(X,\sigma)= T\!\left(Y;\pm c(r_1r_2)^{-1},u_1r_1,u_2r_2,\varepsilon\right)
-\delta_1,
\]
where $\delta_1=1$ if the term corresponding to $D=1$ is counted
on the right, and $\delta_1=0$ otherwise.

By Proposition~\ref{en:prop:sieve}, the main term of
$N_2^\pm(X,\sigma)$ is
\[
 \frac{X}{2^{e+\omega(u_1r_1)}r_1r_2\pi^2}\prod_{2\neq p\mid n}\frac p{p+1}
 =\frac{3X}{\pi^2}\prod_{\substack{p\mid n}}c_p(\sigma_p).
\]
Consequently,
\begin{equation}\label{en:eq:pure-count}
 N_2^\pm(X,\sigma)
 =\frac{3}{\pi^2}X\prod_{p\mid n}c_p(\sigma_p)
 +O\!\left(
 \sqrt{\frac{X}{r_1r_2}}\,
 (\log(2u_1r_1))^{5/4}
 \lambda_1(u_1r_1)\lambda_2(u_2r_2)
 \right).
\end{equation}
The contribution of $\delta_1$ can be absorbed into this error term.
Since $\lambda_1(p)/p^{3/16}\to\sqrt2<3/2$ as $p\to\infty$,
we have $\lambda_1(q)\ll q^{3/16}(3/2)^{\omega(q)}$ for every
positive squarefree integer $q$.
Applying this to \eqref{en:eq:pure-count} proves \eqref{en:eq:main}.
\end{proof}

We now turn to the case of a general local specification
$\sigma=(\sigma_p)_{p\mid n}$.
For each odd prime $p\mid n$, let $\sigma_p^{\mathrm{ur}}$ be the
subset of unramified algebras in $\sigma_p$, and let
$\sigma_p^{\mathrm{rm}}$ be the subset of ramified algebras.
Put
$a_p:=\#\sigma_p^{\mathrm{ur}},
 b_p:=\#\sigma_p^{\mathrm{rm}}$.
Call $\sigma_p$ \emph{pure} if $a_pb_p=0$, and \emph{mixed} otherwise.
Thus $\sigma$ is pure precisely when $\sigma_p$ is pure for every
odd prime $p\mid n$.

We group the odd primes $p\mid n$ according to the numbers
$(a_p,b_p)$ of unramified and ramified algebras in $\sigma_p$,
and define the following five factors of $n$:
\begin{equation}\label{en:eq:local-factors}
 u_1:=\prod_{a_p=1}p,\quad
 u_2:=\prod_{\substack{a_p=2\\b_p=0}}p,\quad
 u_2':=\prod_{\substack{a_p=2\\b_p=1}}p,\quad
 r_1:=\prod_{\substack{a_p=0\\b_p=1}}p,\quad
 r_2:=\prod_{\substack{a_p=0\\b_p=2}}p.
\end{equation}
These five integers are pairwise coprime, and their product equals
the odd part of $n$.
If $\sigma$ is pure, then $u_2'=1$, and the definitions of
$u_1,u_2,r_1,r_2$ agree with those in Section~\ref{en:sec:statement}.
For each odd prime $p\mid n$, set
$c_p(\sigma_p):=(a_p p+b_p)/(2(p+1))$;
the definition of the local density at $p=2$ remains unchanged.

Bhargava, Thorne, and the author~\cite[Proposition 4.2]{en:BTT}
proved, for general local specifications $\sigma$ and every
$\epsilon>0$, that
\begin{equation}\label{en:eq:BTT-general}
 N_2^\pm(X,\sigma)
 =\frac{3}{\pi^2}X\prod_{p\mid n}c_p(\sigma_p)
+ O_{\epsilon}\!\left(
 X^{1/2}\,u_1^{1/4}(u_2u_2')^{1/4}r_1^{-1/4}r_2^{-1/2}
 (u_1u_2u_2'r_1r_2)^{\epsilon}
 \right),
\end{equation}
where the implied constant depends only on $\epsilon$.
The following corollary generalizes Theorem~\ref{en:thm:main} and
improves this estimate.
\begin{corollary}\label{en:cor:general}
For every local specification $\sigma$ and $X>0$, we have
\begin{equation}\label{en:eq:general}
 N_2^\pm(X,\sigma)
 =\frac{3}{\pi^2}X\prod_{p\mid n}c_p(\sigma_p)
 +O\!\left(
 \frac{X^{1/2}\,u_1^{3/16}}{r_1^{5/16}r_2^{1/2}}
 \left(\frac32\right)^{\omega(u_1r_1u_2')}
 (\log(2u_1r_1))^{5/4}\lambda_2(u_2r_2)
 \right).
\end{equation}
\end{corollary}
\begin{proof}
As in the proof of Theorem~\ref{en:thm:main}, it suffices to consider
the case where $2\mid n$ and $\sigma_2$ consists of a single element.
Let $m$ be the product of the odd primes $p$ for which $\sigma_p$ is mixed.
By definition, $u_2'\mid m\mid u_1u_2'$.
For each divisor $t$ of $m$, define a local specification
$\sigma^{(t)}=(\sigma_p^{(t)})_{p\mid n}$
with the same support $n$ as $\sigma$ by setting, for $p\mid m$,
\[
 \sigma_p^{(t)}:=
 \begin{cases}
  \sigma_p^{\mathrm{rm}},&p\mid t,\\
  \sigma_p^{\mathrm{ur}},&p\mid m/t
 \end{cases}
\]
and taking $\sigma_p^{(t)}:=\sigma_p$ for $p\nmid m$.
Splitting each mixed condition $\sigma_p$ into its unramified and
ramified parts yields the disjoint union
$\mathcal F(\sigma)=\sqcup_{t\mid m}\mathcal F(\sigma^{(t)})$.
Therefore,
\begin{equation}\label{en:eq:count-decomposition}
N_2^\pm(X,\sigma)
=\sum_{t\mid m}N_2^\pm(X,\sigma^{(t)}).
\end{equation}
Since each $\sigma^{(t)}$ is pure,
we can apply \eqref{en:eq:pure-count} to each term on the right.
The main terms sum to the desired main term by additivity of the
local densities.

For $\sigma^{(t)}$, let $n_1^{(t)}$ be the product of the odd
primes $p\mid n$ at which exactly one local algebra is allowed,
and let $n_2^{(t)}$ be the corresponding product where exactly two
algebras are allowed.
Their product is the odd part of $n$.
For each $t\mid m$, put $t_0:=(t,u_2')$.
If $p\mid u_2'$, then $\#\sigma_p^{(t)}=1$ precisely when the
ramified part is selected, that is, when $p\mid t$.
This observation and \eqref{en:eq:local-factors} give
$n_1^{(t)}\mid u_1r_1t_0$.
Hence
\[
 \log(2n_1^{(t)})
 \le \log(2u_1r_1t_0)
 \ll \log(2u_1r_1)\cdot t_0^{1/5}.
\]
The product of the odd primes where $\sigma^{(t)}$ prescribes
ramification is $r_1r_2t$.
Summing the errors in \eqref{en:eq:pure-count} for
$N_2^\pm(X,\sigma^{(t)})$ over $t\mid m$, we obtain
\begin{align}
\left|N_2^\pm(X,\sigma)-\frac{3X}{\pi^2}\prod_{p\mid n}c_p(\sigma_p)\right|
&\ll \sqrt{\frac{X}{r_1r_2}}\sum_{t\mid m}t^{-1/2}
 (\log(2n_1^{(t)}))^{5/4}\lambda_1(n_1^{(t)})\lambda_2(n_2^{(t)})\notag\\
 &\ll \sqrt{\frac{X}{r_1r_2}}(\log(2u_1r_1))^{5/4}
 \sum_{t\mid m}\frac{t_0^{1/4}}{t^{1/2}}\lambda_1(n_1^{(t)})\lambda_2(n_2^{(t)})\label{en:eq:mixed-errors-before}
\end{align}
The desired error bound now follows from the next lemma.
\end{proof}

\begin{lemma}\label{en:lem:mixed-sum}
With the notation in the proof of Corollary~\ref{en:cor:general},
we have
\begin{equation}\label{en:eq:absorb}
 \sum_{t\mid m}\frac{t_0^{1/4}}{t^{1/2}}\lambda_1(n_1^{(t)})\lambda_2(n_2^{(t)})
 \ll \left(\frac32\right)^{\omega(u_1r_1u_2')}
 (u_1r_1)^{3/16}\lambda_2(u_2r_2).
\end{equation}
\end{lemma}
\begin{proof}
By multiplicativity, the divisor sum on the left-hand
side of \eqref{en:eq:absorb} factors as
\begin{equation}\label{en:eq:mixed-product}
\prod_{\substack{2\neq p\mid n\\p\nmid m}}\lambda_{a_p+b_p}(p)
 \prod_{p\mid m/u_2'}\bigl(\lambda_{a_p}(p)+p^{-1/2}\lambda_{b_p}(p)\bigr)
 \prod_{p\mid u_2'}\bigl(\lambda_2(p)+p^{-1/4}\lambda_1(p)\bigr).
\end{equation}
Indeed, for $p\mid m$ with $p\nmid u_2'$, the local factor is
the sum of $\lambda_{a_p}(p)$, arising when $p\nmid t$, and
$p^{-1/2}\lambda_{b_p}(p)$, arising when $p\mid t$.
For $p\mid u_2'$, we have $(a_p,b_p)=(2,1)$, and the extra weight
$t_0^{1/4}$ gives the factor
$\lambda_2(p)+p^{-1/4}\lambda_1(p)$.
If $p\nmid m$, only $p\nmid t$ occurs; the factor is a single term,
which can be written as $\lambda_{a_p+b_p}(p)$ because $a_pb_p=0$.

Since $\lambda_2(p)<\kappa(p)$, we have
$\lambda_2(p)<\lambda_1(p)$.
Thus, for $p\mid u_1$, the factor in \eqref{en:eq:mixed-product}
in the mixed case satisfies
\[
\lambda_1(p)+p^{-1/2}\lambda_{b_p}(p)
\leq \lambda_1(p)(1+p^{-1/2})
=\lambda_1(p)\lambda_2(p),
\]
whereas in the pure case it is $\lambda_1(p)$.
In either case, the factor is at most $\lambda_1(p)\lambda_2(p)$.
The factor for $p\mid u_2$ in \eqref{en:eq:mixed-product} is
$\lambda_2(p)$, and that for $p\mid r_i$ is $\lambda_i(p)$.
Consequently,
\begin{align*}
 \sum_{t\mid m}\frac{t_0^{1/4}}{t^{1/2}}\lambda_1(n_1^{(t)})\lambda_2(n_2^{(t)})
& \leq\lambda_1(r_1)\lambda_2(u_2r_2)
 \prod_{p\mid u_1}\lambda_1(p)\lambda_2(p)
 \prod_{p\mid u_2'}\bigl(\lambda_2(p)+p^{-1/4}\lambda_1(p)\bigr)\\
& \leq\lambda_2(u_2r_2)
 \prod_{p\mid u_1r_1}\lambda_1(p)\lambda_2(p)
 \prod_{p\mid u_2'}\bigl(\lambda_2(p)+p^{-1/4}\lambda_1(p)\bigr).
\end{align*}
As $p\to\infty$ through the primes,
\[
\frac{\lambda_1(p)\lambda_2(p)}{p^{3/16}}\longrightarrow\sqrt2<\frac32,
 \qquad
 \lambda_2(p)+p^{-1/4}\lambda_1(p)\longrightarrow1<\frac32.
\]
These bounds prove \eqref{en:eq:absorb}.
\end{proof}
\begin{remark}
The factor $3/2$ in the error terms of \eqref{en:eq:main}
and \eqref{en:eq:general} (and also on the right-hand side of \eqref{en:eq:absorb})
may be replaced by any fixed real number
$B>\sqrt2$ (for example, $B=1.42$).
The implied constant may then depend on $B$.
\end{remark}

\address{Department of Mathematics, Graduate School of Science, Kobe University, Kobe, 657-8501, Japan}
\email{tani@math.kobe-u.ac.jp}


\clearpage
\renewcommand{\bilingualversion}{ja}
\setcounter{section}{0}
\setcounter{subsection}{0}
\setcounter{subsubsection}{0}
\setcounter{equation}{0}
\setcounter{theorem}{0}
\setcounter{footnote}{0}
\renewcommand{\bilingualtheoremname}{定理}
\renewcommand{\bilinguallemmaname}{補題}
\renewcommand{\bilingualpropositionname}{命題}
\renewcommand{\bilingualcorollaryname}{系}
\renewcommand{\bilingualremarkname}{注}
\renewcommand{\proofname}{証明}
\renewcommand{\abstractname}{概要}
\renewcommand{\refname}{参考文献}
\makeatletter


\def\baselinestretch{1}
\renewcommand{\normalsize}{%
  \@setfontsize\normalsize\@xipt{15.5}%
  \abovedisplayskip 11\p@ \@plus3\p@ \@minus6\p@
  \abovedisplayshortskip \z@ \@plus3\p@
  \belowdisplayshortskip 6.5\p@ \@plus3.5\p@ \@minus3\p@
  \belowdisplayskip \abovedisplayskip
  \let\@listi\@listI}
\renewcommand{\small}{%
  \@setfontsize\small\@xpt\@xiipt
  \abovedisplayskip 10\p@ \@plus2\p@ \@minus5\p@
  \abovedisplayshortskip \z@ \@plus3\p@
  \belowdisplayshortskip 6\p@ \@plus3\p@ \@minus3\p@
  \def\@listi{\leftmargin\leftmargini
    \topsep 6\p@ \@plus2\p@ \@minus2\p@
    \parsep 3\p@ \@plus2\p@ \@minus\p@
    \itemsep\parsep}%
  \belowdisplayskip\abovedisplayskip}
\renewcommand{\footnotesize}{%
  \@setfontsize\footnotesize\@ixpt{11}%
  \abovedisplayskip 8\p@ \@plus2\p@ \@minus4\p@
  \abovedisplayshortskip\z@ \@plus\p@
  \belowdisplayshortskip 4\p@ \@plus2\p@ \@minus2\p@
  \def\@listi{\leftmargin\leftmargini
    \topsep 4\p@ \@plus2\p@ \@minus2\p@
    \parsep 2\p@ \@plus\p@ \@minus\p@
    \itemsep\parsep}%
  \belowdisplayskip\abovedisplayskip}
\renewcommand{\scriptsize}{\@setfontsize\scriptsize\@viiipt{9.5}}
\renewcommand{\tiny}{\@setfontsize\tiny\@vipt\@viipt}
\renewcommand{\large}{\@setfontsize\large\@xiipt{17}}
\renewcommand{\Large}{\@setfontsize\Large\@xivpt{21}}
\renewcommand{\LARGE}{\@setfontsize\LARGE\@xviipt{25}}
\renewcommand{\huge}{\@setfontsize\huge\@xxpt{28}}
\renewcommand{\Huge}{\@setfontsize\Huge\@xxvpt{33}}
\normalsize\selectfont

\setlength\parindent{10.53625pt}
\setlength\topskip{8.51457pt}
\setlength\maxdepth{.5\topskip}
\global\@maxdepth\maxdepth
\setlength\headheight{12pt}
\setlength\headsep{25pt}
\setlength\footskip{30pt}
\setlength\lineskip{1pt}
\setlength\normallineskip{1pt}
\setlength\parskip{0pt plus 1pt}
\setlength\footnotesep{7.7pt}
\setlength{\skip\footins}{10pt plus 4pt minus 2pt}
\raggedbottom


\renewcommand{\section}{\@startsection{section}{1}{\z@}%
  {23.25pt plus 7.75pt minus 3.1pt}%
  {7.75pt plus 4.65pt}%
  {\normalfont\Large\bfseries}}
\renewcommand{\subsection}{\@startsection{subsection}{2}{\z@}%
  {23.25pt plus 7.75pt minus 3.1pt}%
  {7.75pt plus 4.65pt}%
  {\normalfont\large\bfseries}}
\renewcommand{\subsubsection}{\@startsection{subsubsection}{3}{\z@}%
  {23.25pt plus 7.75pt minus 3.1pt}%
  {7.75pt plus 4.65pt}%
  {\normalfont\normalsize\bfseries}}

\setlength\leftmargini{2.5em}
\setlength\leftmarginii{2.2em}
\setlength\leftmarginiii{1.87em}
\setlength\leftmarginiv{1.7em}
\setlength\leftmarginv{1em}
\setlength\leftmarginvi{1em}
\setlength\leftmargin{\leftmargini}
\setlength\labelsep{.5em}
\setlength\labelwidth{\leftmargini}
\addtolength\labelwidth{-\labelsep}
\setlength\partopsep{3pt plus 1pt minus 1pt}
\def\@listi{\leftmargin\leftmargini
  \parsep 4.5\p@ \@plus2\p@ \@minus\p@
  \topsep 9\p@ \@plus3\p@ \@minus5\p@
  \itemsep 4.5\p@ \@plus2\p@ \@minus\p@}
\let\@listI\@listi
\@listi
\def\@listii{\leftmargin\leftmarginii
  \labelwidth\leftmarginii\advance\labelwidth-\labelsep
  \topsep 4.5\p@ \@plus2\p@ \@minus\p@
  \parsep 2\p@ \@plus\p@ \@minus\p@
  \itemsep\parsep}
\def\@listiii{\leftmargin\leftmarginiii
  \labelwidth\leftmarginiii\advance\labelwidth-\labelsep
  \topsep 2\p@ \@plus\p@ \@minus\p@
  \parsep\z@
  \partopsep\p@ \@plus\z@ \@minus\p@
  \itemsep\topsep}
\renewenvironment{quotation}
  {\list{}{\listparindent 1.5em%
    \itemindent\listparindent
    \rightmargin\leftmargin
    \parsep\z@ \@plus\p@}%
    \item\relax}
  {\endlist}
\renewenvironment{abstract}
  {\small
    \begin{center}%
      {\bfseries\abstractname\vspace{-.5em}\vspace{\z@}}%
    \end{center}%
    \quotation}
  {\endquotation}

\setlength\bibindent{1.5em}
\normalsize\selectfont
\makeatother
\xeCJKsetup{
  CJKglue={\hspace{0pt plus .4pt minus .5pt}},
  CJKecglue={\hspace{2.40553pt plus 1pt minus 1pt}},
  AllowBreakBetweenPuncts=true,
  CheckSingle=true,
  WidowPenalty=500,
  PunctStyle=fullwidth,
  RubberPunctSkip=false
}
\makeatletter
\newbox\JAFixedGlyphBox
\protected\def\JAFixedMetricGlyph#1{%
  \setbox\JAFixedGlyphBox=\hbox{\XeTeXglyph\XeTeXcharglyph`#1\relax}%
  \ht\JAFixedGlyphBox=.777586\dimexpr\f@size pt\relax
  \dp\JAFixedGlyphBox=.138854\dimexpr\f@size pt\relax
  \box\JAFixedGlyphBox
  \csname xeCJK_CJK_and_Boundary:w\endcsname}
\let\CJKsymbol\JAFixedMetricGlyph
\protected\def\JAFixedPunctGlyph#1{%
  \setbox\JAFixedGlyphBox=\hbox{\XeTeXglyph\XeTeXcharglyph`#1\relax}%
  \ht\JAFixedGlyphBox=.777586\dimexpr\f@size pt\relax
  \dp\JAFixedGlyphBox=.138854\dimexpr\f@size pt\relax
  \box\JAFixedGlyphBox
  \csname __xeCJK_punct_if_right:NTF\endcsname #1
    {\csname xeCJK_FullRight_and_Boundary:\endcsname}
    {\csname xeCJK_FullLeft_and_Boundary:\endcsname}}
\let\CJKpunctsymbol\JAFixedPunctGlyph

\makeatother

\makeatletter
\DeclareRobustCommand\em{%
  \@nomath\em
  \ifdim\fontdimen\@ne\font>\z@
    \CJKfamily{\CJKrmdefault}\upshape
  \else
    \CJKfamily{\CJKsfdefault}\itshape
  \fi}
\makeatother
\DeclareRobustCommand{\S}{\ifmmode\mathsection\else\UseTextSymbol{OMS}{\textsection}\fi}
\title{局所条件を付けた二次体の一様な計数}
\author{谷口隆}
\date{}
\pdfbookmark[0]{日本語原著}{ja:version}
\bilingualmaketitle


\begin{abstract}
局所条件付き二次体計数に関する
Bhargava, Thorne と著者の一様評価を改良する．
Burgess 評価と P\'olya--Vinogradov の不等式を組み合わせ，
誤差項の局所条件を課す素数への依存を改善する．
さらに，従来はそれらの素数の積の
任意の $\epsilon$ 乗に吸収されていた因子を明示し，
暗黙の定数が絶対定数である評価を与える．
\end{abstract}

\section{序}\label{ja:sec:statement}

素数 $p$ に対し，$A_p$を二次エタール $\Q_p$-代数の同型類の集合とする．
$n$ を正の無平方数とする．
各素数 $p\mid n$ に対し，$A_p$ の空でない真部分集合 $\sigma_p$ を一つずつ指定する．
この族 $\sigma=(\sigma_p)_{p\mid n}$ を二次体の\emph{局所条件}と呼び，
$n$ を $\sigma$ の{\emph 台}と呼ぶ．
二次体$F$は，すべての$p\mid n$に対して
$F\otimes \Q_p\in\sigma_p$であるとき，
局所条件$\sigma$をみたすという．
局所条件$\sigma$をみたす二次体の同型類の集合を
$\mathcal F(\sigma)$と書き，
各符号$\pm$と$X>0$ に対して
\[
 N_2^\pm(X,\sigma)
 :=\#\{F\in \mathcal F(\sigma)\mid 0<\pm\Disc(F)<X\}
\]
とおく．
本稿の主定理は、この計数関数の $X$ と $\sigma$ についての一様な評価である．
すべての奇素数 $p\mid n$ に対して，
$\sigma_p$ が不分岐代数と分岐代数の両方を含んではいないとき，
局所条件 $\sigma$ を\emph{純粋}という．
本節では，応用上重要な、
$\sigma$が純粋な場合に主定理を述べる．
一般の局所条件については
第\CJKecglue{}\ref{ja:sec:local}\CJKecglue{}節の系~\ref{ja:cor:general} を参照されたい．

主定理を述べるために，
純粋な局所条件 $\sigma$ に伴う量をいくつか定義する．
$i=1,2$ に対し，$\sigma_p$ が $i$ 個の不分岐代数からなる
奇素数 $p\mid n$ の積を $u_i$，
$i$ 個の分岐代数からなる
奇素数 $p\mid n$ の積を $r_i$ とする．
$n$ が奇数なら $n=u_1u_2r_1r_2$，偶数なら $n=2u_1u_2r_1r_2$ である．
各奇素数 $p\mid n$ における局所密度 $c_p(\sigma_p)$ は，
$p\mid u_i$ のとき $c_p(\sigma_p):=ip/(2(p+1))$，
$p\mid r_i$ のとき $c_p(\sigma_p):=i/(2(p+1))$ と定める（$i=1,2$）．
$2\mid n$ の場合，
判別式の $2$-進付値が $e$ である$\Q_2$上の$2$次エタール代数の
局所密度を$1/(3\cdot 2^e)$とし，
それを$\sigma_2$ の元について足し合わせたものを$c_2(\sigma_2)$とする．

無平方な正整数 $q$ に対し，
$\omega(q)$ で $q$ の相異なる素因子の個数を表し，
$\lambda_2(q):=\prod_{p\mid q}(1+p^{-1/2})$ とおく．
本論文の主結果は次のとおりである．
\begin{theorem}\label{ja:thm:main}
純粋な局所条件 $\sigma$ と $X>0$ に対して次が成り立つ。
\begin{equation}\label{ja:eq:main}
 N_2^\pm(X,\sigma)
 =\frac{3}{\pi^2}X\prod_{p\mid n}c_p(\sigma_p)
 +O\!\left(
 \frac{X^{1/2}\,u_1^{3/16}}{r_1^{5/16}r_2^{1/2}}\,
 \left(\frac32\right)^{\omega(u_1r_1)}
 (\log(2u_1r_1))^{5/4}\lambda_2(u_2r_2)
 \right).
\end{equation}
暗黙の定数は絶対定数である．
\end{theorem}

局所条件を課さない $n=1$ の場合の計数は，
Gegenbauer~\cite{ja:Gegenbauer} による無平方数の計数に帰着される．
固定された局所条件 $\sigma$ に対しては，
合同類内の無平方整数を\cite{ja:Gegenbauer}と同様に数えることができ，
\eqref{ja:eq:main} の誤差項を
$O_\sigma(X^{1/2})$ とした評価が得られる．
この誤差項の $\sigma$ への依存を精密に評価するのが
本稿の目的である。

局所条件への依存を明示した一様評価にも古くから研究がある．
本稿の $u_1=r_1=r_2=1$ の場合に対応する互いに素条件付き計数については，
Cohen~\cite[Lemma 5.2]{ja:Cohen} が一様評価を与え，
Suryanarayana~\cite[Corollary 1]{ja:Suryanarayana} や
Robert と Tenenbaum~\cite[(10.1)]{ja:RT}
が
誤差項の法への依存を精密化した．
定理~\ref{ja:thm:main} は，この場合の誤差項の法への依存をさらに改善する．
また，合同類内の無平方整数の計数については，
CohenとRobinson~\cite[Corollary 1]{ja:CR}
や
Schwarz~\cite[Lemma 8]{ja:Schwarz} が，
法と剰余類に関する一様評価を与えている．
Bhargava, Thorne と著者~\cite[Proposition 4.2]{ja:BTT} は
任意の$\sigma$の場合を扱い、
純粋な$\sigma$に対しては、\eqref{ja:eq:main}\CJKecglue{}の誤差項について
\begin{equation}\label{ja:eq:BTT}
 O_{\epsilon}\!\left(
 {X^{1/2}}\,u_1^{1/4}u_2^{1/4}r_1^{-1/4}r_2^{-1/2}
 (u_1u_2r_1r_2)^{\epsilon}
 \right)
\end{equation}
の評価を得た．
ただし $\epsilon>0$ は任意で，暗黙の定数は $\epsilon$ にのみ依存する．
この評価では，分岐する局所拡大を個別に指定する条件を許しながら，
分岐を指定する素数の積に関する明示的な減衰因子
$r_1^{-1/4}r_2^{-1/2}$ が得られている．
\eqref{ja:eq:main} は，その冪因子
$u_1^{1/4}u_2^{1/4}r_1^{-1/4}$ を
$u_1^{3/16}r_1^{-5/16}$ に改善するとともに，
\eqref{ja:eq:BTT}\CJKecglue{}では
台\ の $\epsilon$ 乗に吸収されていた因子を明示している．

定理\CJKecglue{}\ref{ja:thm:main}\CJKecglue{}の証明の中心となるのは，
無平方な整数上の指標和\CJKecglue{}\eqref{ja:eq:T_h}\CJKecglue{}の評価である．
各素数$p\mid u_1r_1$にお\linebreak[4]ける局所条件の指示関数を二次指標で表し、
積を展開すると，このような和が現れる．
この種の和に\linebreak[4]ついては，
Granville と Soundararajan~\cite[Lemma 4.1]{ja:GS} が
P\'olya--Vinogradov の不等式による評価\linebreak[4]を，
Munsch~\cite[Lemma 2.3]{ja:Munsch} が
Burgess~\cite{ja:Burgess} の不等式による評価を与えている．
本稿では，M\"obius\linebreak[4] 反転の後に現れる\CJKecglue{}
\eqref{ja:eq:T_h-Mobius}\CJKecglue{}の内側の指標和に対し，
和を取る区間の長さに応じて，
P\'olya--Vinogradov の\linebreak[4]不等式と
Jain-Sharma, Khale と Liu~\cite[Theorem 1.1]{ja:JKL}による
Burgess 評価を使い分ける．
また，局\linebreak[4]所条件に依存する因子を最後まで追跡し，
誤差項の局所条件への依存を、素数ごとの因子の積として\linebreak[4]明示する．

Thorne と著者~\cite{ja:TT}，および Bhargava, Thorneと著者~\cite{ja:BTT} による
局所条件付き三次体計数の一様評価は，
数論統計のさまざまな問題に応用されている．
こうした数体計数の証明では，可約エタール代数の寄与を分離するために，
局所条件付きの低次数体計数が必要となることがある．
例えば，\cite[\S4]{ja:BTT} では，可約三次エタール代数
$\Q\times F$（$F/\Q$ は二次体）の寄与を評価するために，
局所条件付き二次体計数が用いられる．
また，Shankar と Tsimerman~\cite[\S3.1]{ja:ST} は，
四次体計数に現れる可約四次エタール代数の寄与を評価する際に，
同様の二次体計数を用いている．
本稿は，これらの議論で用いられる二次体の一様計数評価を改善する．
局所条件付き二次体計数の他の応用例としては，
Pollack~\cite[\S5.2]{ja:Pollack} による
二次体の最小惰性素数の平均値の研究や，
Ellenberg, Pierce と Wood~\cite[\S6]{ja:EPW} による
イデアル類群の $\ell$-等分点の評価が挙げられる．
後者では，二次体計数が篩の入力として用いられている．
本稿の定理~\ref{ja:thm:main} とその一般化である系~\ref{ja:cor:general} も，
数論統計へのさらなる応用をもつことが期待される．

本稿の構成は次のとおりである．
第2節では Burgess 評価と P\'olya--Vinogradov の不等式を組み合わせることで，
指標和の絶対値の総和の一様評価を与える．
第3節では，これと M\"obius 反転を用いて，
合同条件とLegendre記号の値を指定した無平方数の
一様な計数公式を得る．
第4節では，二次体判別式の計数を第3節の計数問題に帰着させ，
定理\CJKecglue{}\ref{ja:thm:main}\CJKecglue{}を証明する．
さらに，一般の局所条件に対する
系~\ref{ja:cor:general} を導く．
特に断らない限り、
Landauの$O$やVinogradovの$\ll$の暗黙定数は絶対定数である。

\subsection*{謝辞}
この研究は科研費JP23K22386, JP24K21512と京都大学に設置された国際共同利用・共同研究拠点である数理解析研究所の助成を受けた。
本研究の一部は、著者がLodha Mathematical Sciences Instituteで開催された
Arithmetic Statisticsをテーマとするプログラムに参加している間に行われた。
著者は、極めて優れた研究環境を提供してくださったプログラムの主催者
およびスタッフの方々に心から感謝申し上げる。

\subsection*{翻訳およびAIの利用に関する注記}
著者はこの論文を母語である日本語でまず執筆し、
その後、2026年10月から11月にかけて、
ChatGPT（GPT-6）の支援を受けて英語版へ翻訳した。
また、著者はこの間、数学的な記法や用語についてもChatGPTに相談した。
著者は翻訳全体を慎重に確認・修正しており、論文の内容について全責任を負う。
日本語版と英語版はいずれも著者のウェブサイトから取得できる。

\section{指標和の評価}\label{ja:sec:characters}

本節では，$\chi$ を非自明な原始 Dirichlet 指標とし，その導手を $q$ とする．
非負整数 $N$ と実数 $b$ に対し，$N$ 項の指標和
\[
 S(\chi,N;b):=\sum_{b<j\le b+N}\chi(j)
\]
を定める．
正整数$q$に対して
$\omega(q)$を$q$の相異なる素因数の個数,
$\tau(q)$ を $q$ の正の約数の個数，
$\varphi(q)$ を $(\Z/q\Z)^\times$の位数とし，
また
\begin{equation}\label{ja:eq:kappa}
 \kappa(q):=q^{3/16}(2^{\omega(q)}\tau(q))^{3/4}
 \left(\frac q{\varphi(q)}\right)^{1/2}
\end{equation}
と定める．
$2^{\omega},\tau,\varphi$はいずれも乗法的関数なので，
$\kappa$も乗法的関数である。
Burgess~\cite{ja:Burgess} の指標和評価の明示的な形として，
Jain-Sharma, Khale と Liu~\cite[Theorem 1.1]{ja:JKL}
による次の評価を用いる．
\begin{lemma}[Burgess 評価 \cite{ja:JKL}]\label{ja:lem:all-q}
$N\le q^{5/8}$ のとき，次が成り立つ．
\begin{equation}\label{ja:eq:JKL-usable}
 |S(\chi,N;b)|\ll\sqrt N\,\kappa(q)(\log q)^{1/4}.
\end{equation}
\end{lemma}
\cite[Theorem 1.1]{ja:JKL} における導手が十分大きいという仮定は，
暗黙の定数を大きくすることで除ける．
次はよく知られたP\'olya--Vinogradovの不等式である．
証明はたとえば Davenport~\cite[\S23]{ja:Davenport} を参照のこと．
\begin{lemma}[P\'olya--Vinogradov]\label{ja:lem:PV}
任意の$\chi,N,b$に対し，次が成り立つ．
\[
 |S(\chi,N;b)|\ll\sqrt q\,\log q.
\]
\end{lemma}

Burgess 評価と P\'olya--Vinogradov の不等式を組み合わせて，次を示す．
\begin{proposition}\label{ja:prop:aggregate}
$Y>0$ とし，
各整数 $a\ge1$ に対し、
非負整数$N_a$ が
$N_a\le Y/a^2$ を満たしているとする．このとき
\begin{equation}\label{ja:eq:aggregate}
 \sum_{a\ge1}|S(\chi,N_a;b_a)|
 \ll\sqrt Y\,\kappa(q)(\log q)^{5/4}
\end{equation}
である．ただし各$b_a$は任意の実数である。
\end{proposition}
\begin{proof}
$
 A=\sqrt Y\,q^{-5/16}
$
とおく．$a\leq A$ の部分には補題~\ref{ja:lem:PV} を用いて
\[
 \sum_{a\leq A}|S(\chi,N_a;b_a)|
 \ll A\sqrt q\,\log q
 =\sqrt Y\,q^{3/16}\log q
 \ll\sqrt Y\,\kappa(q)(\log q)^{5/4}
\]
\begingroup
\xeCJKsetup{CJKglue={\hspace{0pt plus .4pt minus .7pt}}}
を得る．
$a> A$とすると$N_a\leq Y/A^2=q^{5/8}$である。
また、$N_a>0$ なら，整数性より$N_a\geq1$だから\linebreak[4]
$a\le\sqrt Y$ である．従って補題~\ref{ja:lem:all-q} より\endgroup{}
\[
 \sum_{a> A}|S(\chi,N_a;b_a)|
 \ll     \sum_{A< a\le\sqrt Y}
 \sqrt{N_a}\,\kappa(q)(\log q)^{1/4}
 \leq\sqrt Y\,\kappa(q)(\log q)^{1/4}
      \sum_{A< a\le\sqrt Y}\frac1a.
\]
\par\removelastskip\newpage\noindent
最後の和は
\[
 \sum_{A< a\le\sqrt Y}\frac1a
 \le1+\log\frac{\sqrt Y}{A}
 =1+\frac5{16}\log q
\]
で抑えられる．これを Burgess 側の評価に代入し，
P\'olya--Vinogradov 側の評価と合わせればよい．
\end{proof}

\section{平方因子を除く基本計数}\label{ja:sec:sieve}

$n_1,n_2$ を互いに素な正の奇数の無平方整数，
$c\in(\Z/8\Z)^\times$ とする．
$\varepsilon=(\varepsilon(p))_{p\mid n_1}$
は各$p\mid n_1$における符号
$\varepsilon(p)\in\{+1,-1\}$
の族とする．
$Y>0$に対し，$T(Y;c,n_1,n_2,\varepsilon)$ で，
次の条件を満たす整数 $0<d<Y$ の個数を表す．
\begin{center}
$d$は無平方数;
\quad
$d\equiv c\pmod8$;
\quad
$(d,n_1n_2)=1$;
\quad
すべての$p\mid n_1$に対し
$\leg dp=\varepsilon(p)$.
\end{center}
定義から、この計数関数は次の表示をもつ．
\begin{equation}\label{ja:eq:T-definition}
 T(Y;c,n_1,n_2,\varepsilon)=
 \sum_{\substack{0<d<Y\\ d\equiv c\ (8)\\ (d,n_1n_2)=1}}
 \mu(d)^2\prod_{p\mid n_1}\frac{1+\varepsilon(p)\leg dp}{2}
\end{equation}
ここで$\mu$ は M\"obius 関数である．
無平方な正整数上の乗法的関数 $\lambda_1$ を，素数 $p$ における値を
\begin{equation}\label{ja:eq:lambda-def}
\lambda_1(p):=\frac{\kappa(p)+\lambda_2(p)}2
 =\frac{\sqrt2}{(1-p^{-1})^{1/2}}\,p^{3/16}+\frac{1+p^{-1/2}}2
\end{equation}
とすることで定める．本節の目標は次の一様評価である。
\begin{proposition}\label{ja:prop:sieve}
\begin{equation}\label{ja:eq:sieve}
 T(Y;c,n_1,n_2,\varepsilon)=\frac{Y}{2^{\omega(n_1)}\pi^2}\prod_{p\mid n_1n_2}\frac p{p+1}
 +O\!\left(\sqrt Y\,(\log(2n_1))^{5/4}\lambda_1(n_1)\lambda_2(n_2)\right).
\end{equation}
\end{proposition}

\begin{proof}
\begingroup
\xeCJKsetup{CJKglue={\hspace{0pt plus .4pt minus .7pt}}}
各 $h\mid n_1$ に対して
$\varepsilon(h):=\prod_{p\mid h}\varepsilon(p)$ と定めることで，
$\varepsilon$を
$n_1$ の約数上の乗法的関数とみなす．\linebreak[4]
\eqref{ja:eq:T-definition} の積は\endgroup{}
\[
 \prod_{p\mid n_1}\frac{1+\varepsilon(p)\leg dp}{2}
 =\frac{1}{2^{\omega(n_1)}}\sum_{h\mid n_1}\varepsilon(h)\leg dh
\]
と展開できる．
ここで $\leg dh$ は Jacobi 記号である．
したがって，各 $h\mid n_1$ に対して
\begin{equation}\label{ja:eq:T_h}
 T_h(Y):=
 \sum_{\substack{0<d<Y\\d\equiv c\ (8)\\ (d,n_1n_2)=1}}
 \mu(d)^2\leg dh
\end{equation}
とおくと，
\begin{equation}\label{ja:eq:sieve-expand}
 T(Y;c,n_1,n_2,\varepsilon)
 =\frac{1}{2^{\omega(n_1)}}\sum_{h\mid n_1}\varepsilon(h)T_h(Y)
\end{equation}
である．$T_h(Y)$ の $c,n_1,n_2$ への依存は記号から省略する．
次に，各 $T_h(Y)$ に
$\mu(d)^2=\sum_{a^2\mid d}\mu(a)$ を代入し，$d=a^2k$ と置く．
もとの条件から $(a,2n_1n_2)=(k,n_1n_2)=1$ であり，
$a^2\equiv1\pmod8$，$\leg{a^2}h=1$ より，
\[
 T_h(Y)=
 \sum_{\substack{a\le\sqrt Y\\ (a,2n_1n_2)=1}}\mu(a)
   \sum_{\substack{0<k<Y/a^2\\ k\equiv c\ (8)\\ (k,n_1n_2)=1}}\leg kh
\]
となる．
$\leg kh$ は $(k,h)>1$ なら零になるため，内側の和の条件 $(k,n_1n_2)=1$ を
$(k,n_1n_2/h)=1$ に置き換えてよい．
この条件を包除公式
\[
 \mathbf 1_{(k,n_1n_2/h)=1}
 =\sum_{e\mid (n_1n_2/h,k)}\mu(e)
 =\sum_{\substack{e\mid n_1n_2/h\\e\mid k}}\mu(e)
\]
で展開し，$k=e\ell$ と置く．$e$ は奇数なので$(\Z/8\Z)^\times$
で逆元 $e^{-1}$ を持ち，
$k\equiv c\pmod8$ は $\ell\equiv ce^{-1}\pmod8$ となる．
$e$ と $a$ の和の順序を入れ替え，以下を得る：
\begin{equation}\label{ja:eq:T_h-Mobius}
 T_h(Y)=\sum_{e\mid n_1n_2/h}\mu(e)\leg eh
   \sum_{\substack{a\le\sqrt{Y/e}\\ (a,2n_1n_2)=1}}\mu(a)
   \sum_{\substack{0<\ell<Y/ea^2\\\ell\equiv ce^{-1}\ (8)}}\leg\ell h
\end{equation}

まず $h=1$の項$T_1(Y)$ を評価する．
この場合，\eqref{ja:eq:T_h-Mobius} の $\ell$ に関する和は $Y/(8ea^2)+O(1)$ である．
これを$a$ について足すと，$O(1)$ の総和も，
主項の $a$についての和を無限和に延長する誤差も $O(\sqrt{Y/e})$ である．従って，
$\sum_{e\mid n_1n_2}e^{-1/2}=\lambda_2(n_1n_2)$ より，
\begin{align*}
T_1(Y)&=\frac{Y}{8}
   \sum_{e\mid n_1n_2}\frac{\mu(e)}e
   \sum_{\substack{a\ge1\\ (a,2n_1n_2)=1}}\frac{\mu(a)}{a^2}
   +O\!\left(\sqrt Y\,\lambda_2(n_1n_2)\right)
 \\
 &=
 \frac{Y}{8}
 \prod_{p\mid n_1n_2}(1-p^{-1})\prod_{p\nmid2n_1n_2}(1-p^{-2})
 +O\!\left(\sqrt Y\,\lambda_2(n_1n_2)\right).
\end{align*}
$\prod_{p\ne2}(1-p^{-2})=8/\pi^2$ より，
$2^{-\omega(n_1)}T_1(Y)$ の主項は \eqref{ja:eq:sieve} の主項と一致する．

次に $h>1$ の項$T_h(Y)$を評価する．
$e\mid n_1n_2/h$に対し、
$r\equiv ce^{-1}\pmod8$となる$r\in\Z$ を選び、
さらに$8b\equiv r\pmod h$ を満たす整数 $b$ を選ぶ。
$N_a$ を $0<\ell<Y/ea^2$，$\ell\equiv r\pmod8$ を満たす整数 $\ell$ の個数とする．
$\ell=r+8j$ と書けば，$\ell>0$ は $j>-r/8$ と同値であり，
$\leg\ell h=\leg8h\leg{j+b}h$ である．従って
\[
 \sum_{\substack{0<\ell<Y/ea^2\\\ell\equiv ce^{-1}\ (8)}}\leg\ell h
 =\leg8h\,S\!\left(\leg{\cdot}h,N_a;b-\frac r8\right),
 \qquad N_a\le\frac{Y}{ea^2}
\]
である。ヤコビ記号$\leg \cdot h$は導手$h$の原始指標なので、
命題~\ref{ja:prop:aggregate} から	
\[
 \sum_{\substack{a\le\sqrt{Y/e}\\ (a,2n_1n_2)=1}}
 \left|S\!\left(\leg{\cdot}h,N_a;b-\frac r8\right)\right|
 \ll\sqrt{\frac Ye}\,\kappa(h)(\log h)^{5/4}
\]
を得る．
よって，$h<2n_1$と
$\sum_{e\mid n_1n_2/h}e^{-1/2}=\lambda_2(n_1n_2/h)$ を用いると，
\begin{equation}\label{ja:eq:nonprincipal}
\sum_{\substack{1\neq h\mid n_1}}|T_h(Y)|
 \ll \sqrt Y\,(\log(2n_1))^{5/4}
     \sum_{\substack{1\neq h\mid n_1}}\kappa(h)\lambda_2(n_1n_2/h).
\end{equation}

$\kappa(1)=1$ なので，$T_1(Y)$ の誤差評価と\CJKecglue{}
\eqref{ja:eq:nonprincipal} を $h\mid n_1$ 全体の和にまとめられる．
$\lambda_2,\kappa$ の乗法性と \eqref{ja:eq:lambda-def} により
\begin{equation}\label{ja:eq:sieve-product}
 \sum_{h\mid n_1}\kappa(h)\lambda_2(n_1n_2/h)
 =\lambda_2(n_2)
 \prod_{p\mid n_1}\bigl(\lambda_2(p)+\kappa(p)\bigr)
 =2^{\omega(n_1)}\lambda_1(n_1)\lambda_2(n_2).
\end{equation}
\eqref{ja:eq:sieve-expand} と合わせれば，\eqref{ja:eq:sieve} が得られる．
\end{proof}

\section{主定理の証明と一般の局所条件}\label{ja:sec:local}
本節ではまず局所条件$\sigma$が純粋な場合の定理\CJKecglue{}\ref{ja:thm:main}\CJKecglue{}を示し、
次に一般の$\sigma$の場合の系\CJKecglue{}\ref{ja:cor:general}\CJKecglue{}を導く。
\begin{proof}[定理~\ref{ja:thm:main} の証明]
まず $p=2$ における局所条件を、
二次体$F$の判別式
$D=\Disc(F)$の合同条件で表す．
$F$に対し、
\[
 e=\operatorname{ord}_2(D)\in\{0,2,3\},\qquad
 c\equiv D/2^e\pmod8,\quad c\in(\Z/8\Z)^\times
\]
とおくと、$F\otimes\Q_2$の同型類は
$(e,c)$のみで決まる．
可能な $(e,c)$ の組は次の $8$ 組である：
\[
(0,1),(0,5),(2,3),(2,7),(3,1),(3,3),(3,5),(3,7)
\]
それぞれに対応する$2$次$\Q_2$代数の密度は
$1/(3\cdot2^e)$ である．

$2\mid n$ かつ $\sigma_2$ が一つの元からなる場合に定理を示せば十分である．
実際，$F\otimes\Q_2$ の同型類によって分割すれば
（必要に応じて台を $n$ から $2n$ に変更して），
一般の場合の計数はこの形の高々 $8$ 個\linebreak[4]の計数の和となり，
$u_1,u_2,r_1,r_2$ は変わらない．
主項は局所密度の加法性（$A_2$ 上の密度の総和は $1$）\linebreak[4]により足し合わされ，
誤差評価は絶対定数倍を除いてそのまま保たれる．

以下，この場合を考え，$\sigma_2$ の唯一の元に
対応する組を $(e,c)$ とする．
このとき $c_2(\sigma_2)=1/(3\cdot 2^e)$ である．

分岐を指定する奇素数の積は $r_1r_2$ である．
したがって、
$N_2^\pm(X,\sigma)$によって数えられる2次体$F$の判別式$D$は、
ある整数$d$によって
$D=\pm 2^e r_1r_2d$
と表される。$D$がみたす条件を$d$の条件で表す。
$Y=X/(2^e r_1r_2)$ とおけば，
$d$は以下をみたす：
\[
\text{$d$は無平方};
\qquad
0<d<Y;
\qquad
(d,n)=1;
\qquad
d\equiv\pm c(r_1r_2)^{-1}\pmod8.
\]
また，奇素数 $p\mid n$ のうち $\#\sigma_p=1$ となるものの積は $u_1r_1$，
$\#\sigma_p=2$ となるものの積は $u_2r_2$ である．
$p\mid u_1r_1$ では条件
$\leg dp=\varepsilon(p)$ が課される。
実際，$p\mid u_1$では $\leg Dp$，
$p\mid r_1$では$\leg{D/p}p$ によって$F\otimes\Q_p$が決まり，
$D=\pm 2^e r_1r_2d$を代入すれば、
これらはいずれも $\leg dp=\pm1$の形の条件になる．
また，$p\mid u_2r_2$ では，すでに課した $(d,p)=1$ のほかに条件はない．
逆に，これらの条件を満たす $d$ は，$D=1$ の場合を除き，
数える二次体 $F$ と一対一に対応する．
したがって，
\[
N_2^\pm(X,\sigma)= T\!\left(Y;\pm c(r_1r_2)^{-1},u_1r_1,u_2r_2,\varepsilon\right)
-\delta_1
\]
である．ただし$\delta_1\in\{0,1\}$で，
$\delta_1=1$は$D=1$に対応する$d$が数えられたときの
補正である．

命題~\ref{ja:prop:sieve} より、$N_2^\pm(X,\sigma)$の主項は
\[
 \frac{X}{2^{e+\omega(u_1r_1)}r_1r_2\pi^2}\prod_{2\neq p\mid n}\frac p{p+1}
 =\frac{3X}{\pi^2}\prod_{\substack{p\mid n}}c_p(\sigma_p)
\]
である．従って，
\begin{equation}\label{ja:eq:pure-count}
 N_2^\pm(X,\sigma)
 =\frac{3}{\pi^2}X\prod_{p\mid n}c_p(\sigma_p)
 +O\!\left(
 \sqrt{\frac{X}{r_1r_2}}\,
 (\log(2u_1r_1))^{5/4}
 \lambda_1(u_1r_1)\lambda_2(u_2r_2)
 \right)
\end{equation}
を得る．$\delta_1$ の影響はこの誤差に吸収できる．
素数 $p\to\infty$ に対して
$\lambda_1(p)/p^{3/16}\to\sqrt2<3/2$ なので，
無平方な正整数 $q$ に対し
$\lambda_1(q)\ll q^{3/16}(3/2)^{\omega(q)}$ である．
これを \eqref{ja:eq:pure-count} に用いれば \eqref{ja:eq:main} を得る．
\end{proof}

次に，一般の局所条件 $\sigma=(\sigma_p)_{p\mid n}$ に対する計数を考える．
奇素数 $p\mid n$ に対し，$\sigma_p$ の元のうち不分岐なもの全体を
$\sigma_p^{\mathrm{ur}}$，分岐するもの全体を $\sigma_p^{\mathrm{rm}}$ と書き，
それぞれの元の個数を
$a_p:=\#\sigma_p^{\mathrm{ur}},
 b_p:=\#\sigma_p^{\mathrm{rm}}
$
とする．$a_pb_p=0$ のとき $\sigma_p$ を純粋，
そうでないとき混合と呼ぶ．
$\sigma$ が純粋であることは、
すべての奇素数 $p\mid n$ で $\sigma_p$ が純粋であることと同値である。

各 $\sigma_p$ に含まれる不分岐代数と分岐代数の個数 $(a_p,b_p)$ に応じて，
奇素数 $p\mid n$ を分類し，$n$ の次の五つの因子を定める：
\begin{equation}\label{ja:eq:local-factors}
 u_1:=\prod_{a_p=1}p,\quad
 u_2:=\prod_{\substack{a_p=2\\b_p=0}}p,\quad
 u_2':=\prod_{\substack{a_p=2\\b_p=1}}p,\quad
 r_1:=\prod_{\substack{a_p=0\\b_p=1}}p,\quad
 r_2:=\prod_{\substack{a_p=0\\b_p=2}}p.
\end{equation}
これら五つの整数は互いに素で，その積は $n$ の奇数部分に等しい．
また，$\sigma$ が純粋な場合には $u_2'=1$ であり，
$u_1,u_2,r_1,r_2$ の定義は第\CJKecglue{}\ref{ja:sec:statement}\CJKecglue{}節のものと一致する．
各奇素数 $p\mid n$ に対し，
局所密度を $c_p(\sigma_p):=(a_p p+b_p)/(2(p+1))$ と定め，
$p=2$ における局所密度の定義は変えない．

Bhargava, Thorne と著者 ~\cite[Proposition 4.2]{ja:BTT} は
一般の局所条件$\sigma$と任意の$\epsilon>0$に対し，
\begin{equation}\label{ja:eq:BTT-general}
 N_2^\pm(X,\sigma)
 =\frac{3}{\pi^2}X\prod_{p\mid n}c_p(\sigma_p)
+ O_{\epsilon}\!\left(
 {X^{1/2}}\,u_1^{1/4}(u_2u_2')^{1/4}r_1^{-1/4}r_2^{-1/2}
 (u_1u_2u_2'r_1r_2)^{\epsilon}
 \right)
\end{equation}
を得た．ただし暗黙定数は$\epsilon$にのみ依存する．
次の系は定理~\ref{ja:thm:main} の一般化であり，
またこの評価を改善する．
\begin{corollary}\label{ja:cor:general}
すべての局所条件$\sigma$と
$X>0$ に対し、以下が成り立つ。
\begin{equation}\label{ja:eq:general}
 N_2^\pm(X,\sigma)
 =\frac{3}{\pi^2}X\prod_{p\mid n}c_p(\sigma_p)
 +O\!\left(
 \frac{{X^{1/2}}\,u_1^{3/16}}{r_1^{5/16}r_2^{1/2}}
 \left(\frac32\right)^{\omega(u_1r_1u_2')}
 (\log(2u_1r_1))^{5/4}\lambda_2(u_2r_2)
 \right).
\end{equation}
\end{corollary}
\begin{proof}
\begingroup
\xeCJKsetup{CJKglue={\hspace{0pt plus .4pt minus .7pt}}}
定理~\ref{ja:thm:main} の証明と同様に，
$2\mid n$ かつ $\sigma_2$ が一つの元からなる場合を示せば十分である．
$\sigma_p$ が混合である奇素数 $p$ の積を $m$ とする．
定義より，$u_2'\mid m\mid u_1u_2'$ である．
$m$ の各約数 $t$ に対して，$\sigma$ と\linebreak[4]同じ台 $n$ をもつ局所条件
$\sigma^{(t)}=(\sigma_p^{(t)})_{p\mid n}$ を，
$p\mid m$においては\endgroup{}
\[
 \sigma_p^{(t)}:=
 \begin{cases}
  \sigma_p^{\mathrm{rm}},&p\mid t,\\
  \sigma_p^{\mathrm{ur}},&p\mid m/t
 \end{cases}
\]
\begingroup
\xeCJKsetup{CJKglue={\hspace{0pt plus .4pt minus .7pt}}}
とし，$p\nmid m$においては
$\sigma_p^{(t)}:=\sigma_p$とすることによって定める．
混合な局所条件$\sigma_p$を
$\sigma_p^{\mathrm{ur}}$と$\sigma_p^{\mathrm{rm}}$に分\linebreak[4]けることにより，
$\mathcal F(\sigma)$は$\mathcal F(\sigma)=\sqcup_{t\mid m}\mathcal F(\sigma^{(t)})$
と無縁和に分割されるから、\endgroup{}
\begin{equation}\label{ja:eq:count-decomposition}
N_2^\pm(X,\sigma)
=\sum_{t\mid m}N_2^\pm(X,\sigma^{(t)})
\end{equation}
である．
各 $\sigma^{(t)}$ は純粋であるから，
主定理の証明で得た \eqref{ja:eq:pure-count} を右辺の各項に適用する．
それぞれの主項を足すと，局所密度の加法性から所望の主項を得る．

局所条件 $\sigma^{(t)}$ のもとで，
許される局所代数が一つである奇素数 $p\mid n$ の積を $n_1^{(t)}$，
二つである奇素数 $p\mid n$ の積を $n_2^{(t)}$ とする．
両者の積は $n$ の奇数部分に等しい。
各 $t\mid m$ に対して $t_0:=(t,u_2')$ とおく．
$p\mid u_2'$ では，分岐側を選ぶ $p\mid t$ の場合に限って
$\#\sigma_p^{(t)}=1$ となる．このことと \eqref{ja:eq:local-factors} より，
$n_1^{(t)}\mid u_1r_1t_0$ である．
したがって，

\[
 \log(2n_1^{(t)})
 \le \log(2u_1r_1t_0)
 \ll \log(2u_1r_1)\cdot t_0^{1/5}
\]
が成り立つ．
$\sigma^{(t)}$ で分岐を指定する奇素数の積は $r_1r_2t$ なので，
$N_2^\pm(X,\sigma^{(t)})$ の \eqref{ja:eq:pure-count} による誤差を
$t\mid m$ について足すことで，
\begin{align}
\left|N_2^\pm(X,\sigma)-\frac{3X}{\pi^2}\prod_{p\mid n}c_p(\sigma_p)\right|
&\ll\sqrt{\frac{X}{r_1r_2}}\sum_{t\mid m}t^{-1/2}
 (\log(2n_1^{(t)}))^{5/4}\lambda_1(n_1^{(t)})\lambda_2(n_2^{(t)})\notag\\
 &\ll\sqrt{\frac{X}{r_1r_2}}(\log(2u_1r_1))^{5/4}
 \sum_{t\mid m}\frac{t_0^{1/4}}{t^{1/2}}\lambda_1(n_1^{(t)})\lambda_2(n_2^{(t)})\label{ja:eq:mixed-errors-before}
\end{align}
を得る．
以下の補題から所望の誤差評価が従う．
\end{proof}

\begin{lemma}\label{ja:lem:mixed-sum}
系~\ref{ja:cor:general} の証明中の記号を用いると，次が成り立つ．
\begin{equation}\label{ja:eq:absorb}
 \sum_{t\mid m}\frac{t_0^{1/4}}{t^{1/2}}\lambda_1(n_1^{(t)})\lambda_2(n_2^{(t)})
 \ll \left(\frac32\right)^{\omega(u_1r_1u_2')}
 (u_1r_1)^{3/16}\lambda_2(u_2r_2)
\end{equation}
\end{lemma}
\begin{proof}
乗法性により，\eqref{ja:eq:absorb} の左辺の約数和は
次のように素数ごとの積に分解される：
\begin{equation}\label{ja:eq:mixed-product}
\prod_{\substack{2\neq p\mid n\\p\nmid m}}\lambda_{a_p+b_p}(p)
 \prod_{p\mid m/u_2'}\bigl(\lambda_{a_p}(p)+p^{-1/2}\lambda_{b_p}(p)\bigr)
 \prod_{p\mid u_2'}\bigl(\lambda_2(p)+p^{-1/4}\lambda_1(p)\bigr)
\end{equation}
実際、$p\mid m$ かつ $p\nmid u_2'$ における因子は，
$p\nmid t$ の場合の $\lambda_{a_p}(p)$ と
$p\mid t$ の場合の $p^{-1/2}\lambda_{b_p}(p)$ の和である．
$p\mid u_2'$ では，$(a_p,b_p)=(2,1)$ と追加の重み $t_0^{1/4}$ により，
因子は $\lambda_2(p)+p^{-1/4}\lambda_1(p)$ となる．
一方，$p\nmid m$ では $p\nmid t$ に限られるので因子は単項であり，
$a_pb_p=0$に注意するとそれは
$\lambda_{a_p+b_p}(p)$ と書ける．

$\lambda_2(p)<\kappa(p)$ より $\lambda_2(p)<\lambda_1(p)$ なので，
\eqref{ja:eq:mixed-product}\CJKecglue{}の$p\mid u_1$ における因子は，混合の場合には
\[
\lambda_1(p)+p^{-1/2}\lambda_{b_p}(p)
\leq \lambda_1(p)(1+p^{-1/2})
=\lambda_1(p)\lambda_2(p)
\]
であり，純粋の場合には $\lambda_1(p)$ であるから、
いずれの場合も$\lambda_1(p)\lambda_2(p)$以下である。
また、\eqref{ja:eq:mixed-product}\CJKecglue{}の$p\mid u_2$における因子は$\lambda_2(p)$であり、
$p\mid r_i$ における因子は $\lambda_i(p)$ である。
従って，
\begin{align*}
 \sum_{t\mid m}\frac{t_0^{1/4}}{t^{1/2}}\lambda_1(n_1^{(t)})\lambda_2(n_2^{(t)})
& \leq\lambda_1(r_1)\lambda_2(u_2r_2)
 \prod_{p\mid u_1}\lambda_1(p)\lambda_2(p)
 \prod_{p\mid u_2'}\bigl(\lambda_2(p)+p^{-1/4}\lambda_1(p)\bigr)\\
& \leq\lambda_2(u_2r_2)
 \prod_{p\mid u_1r_1}\lambda_1(p)\lambda_2(p)
 \prod_{p\mid u_2'}\bigl(\lambda_2(p)+p^{-1/4}\lambda_1(p)\bigr)
\end{align*}
である．素数 $p\to\infty$ において
\[
\frac{\lambda_1(p)\lambda_2(p)}{p^{3/16}}\longrightarrow\sqrt2<\frac32,
 \qquad
 \lambda_2(p)+p^{-1/4}\lambda_1(p)\longrightarrow1<\frac32
\]
である．以上より \eqref{ja:eq:absorb} が従う．
\end{proof}
\begin{remark}
\eqref{ja:eq:main} と \eqref{ja:eq:general} の誤差項，
および \eqref{ja:eq:absorb} の右辺に現れる $3/2$ は，
任意の固定された実数 $B>\sqrt2$（例えば $B=1.42$）に置き換えられる．
このとき暗黙定数は $B$ に依存しうる．
\end{remark}

\par\bigskip
{\small
\noindent 657-8501\ \  兵庫県神戸市灘区六甲台町1-1
\\神戸大学大学院理学研究科数学専攻\par
\noindent
\href{mailto:tani@math.kobe-u.ac.jp}
{\texttt{tani@math.kobe-u.ac.jp}}\par
}


\end{document}